\documentclass[]{article}
\usepackage{url}
\usepackage{enumerate}
\usepackage{graphicx}
\usepackage{array}
\usepackage{amsmath}
\usepackage{epsfig}
\usepackage{amsxtra}
\usepackage{amsbsy}
\usepackage{amsopn}
\usepackage{amsfonts}
\usepackage{amssymb}
\usepackage{mathrsfs}
\usepackage{epsf}
\usepackage{subfigure}
\usepackage{psfrag}
\usepackage{setspace}
\usepackage{xspace}
\usepackage{multirow}
\usepackage{algorithmic}
\usepackage{epstopdf}
\usepackage{color}
\usepackage{subfigure}

\begin{document}
\thispagestyle{empty}

\onehalfspacing

\begin{center} \textbf{Modulo periodic Poisson stable solutions of quasilinear  differential equations}

$^1$Akhmet M., $^2$Tleubergenova M., $^{1,2}$Zhamanshin A.

$^1$Department of Mathematics, Middle East Technical University, 06800, Ankara, Turkey \\
$^2$Department of Mathematics, Aktobe Regional  University, 030000, Aktobe, Kazakhstan
\end{center}

\textbf{Abstract} We introduce  a new type  of recurrence in the  space of continuous  and bounded functions.  The property  is easily  verifiable,  and can  be considered  for   differential  equations.  This time,  the existence and asymptotic stability  of modulo  periodic Poisson stable solutions for quasilinear systems are  proved.  The significant  novelty  of the research  is the numerical simulation  of the functions,  which  is stemmed  from  dynamical traditions of trigonometric functions,  and contributes to applications of Poisson stable oscillations.

\section{Introduction}

 The  theory of differential equations  is,  first of all,    a doctrine on  oscillations  and recurrence,   which  are  basic in  science and technique.   Oscillations   are most preferable  in  engineering \cite{Minorsky},  while recurrence originates in celestial mechanics \cite{Poin}. The ultimate recurrence is the Poisson  stability \cite{Birkhoff,Nemytskii,sell}.  Nowadays, needs for  functions  with irregular behavior are exceptionally   strong in neuroscience and  celestial dynamics,  which  is still  in  the developing  mode.  In  the present  research  we have  decided to combine  periodic dynamics with the  phenomenon  of Poisson  stability.  That is,  one of the simplest form of oscillations is amalgamated with the  most sophisticated recurrence.  We hope  that  the choice  can give a new push  for the nonlinear analysis,  which  faces challenging  problems of the real world and industry.    The present product of the design   are  modulo  periodic Poisson stable stable functions.
 
 In   paper \cite{ak8},  to strengthen  the role of  recurrence as a chaotic ingredient  we have   extended  the Poisson stability  to the unpredictability property.   Thus,  the Poincar\'{e} chaos has been  determined,    and  one can say  that  the \textit{unpredictability implies chaos} now. The unpredictable point  in the  functional  space of  the Bebutov  dynamics is the  unpredictable function \cite{Akhmet1,a4,AkhmetBook,AkhmetExistence,TJM,Akhmet2020,Akhmet2024,iop,Akhmet2025}.   Accordingly,   we provided a dynamical  method,  how to construct  Poisson stable  functions.   Deterministic and stochastic  dynamics  have been  utilized.   Deterministically  unpredictable functions have been  constructed  as solutions of hybrid systems,   consisting of discrete  and differential equations \cite{Akhmet1,AkhmetBook,Akhmet2020}, and  randomly  they  are  results of the Bernoulli process inserted   a linear  differential equation \cite{iop,AkhmetStrings,Akhmet2025}.  Unpredictable oscillations  in  neural networks  have been  researched  in  \cite{AkhmetBook,Akhmet6,Akhmet2023,Akhmet2027,iop}.

 In  papers \cite{Akhmet1,Akhmet2020,Akhmet2024,Akhmet2025} and books \cite{AkhmetBook,iop} discussing  existence of unpredictable solutions,  we have developed  a new method how to  approve Poisson stable solutions,  since unpredictable functions  are a  subset  of Poisson stable functions,  and to  verify  the unpredictability one has to check,   if the Poisson  stability  is valid.  The method  is distinctly different   than  the\textit{ comparability method by character of recurrence} introduced in  \cite{Shcherbakov2} and later has been  realized in  several  articles \cite{Cheban1,Cheban2,Shcherbakov,Shcherbakov1,Shcherbakov3,Shcherbakov4,Shcherbakov5}.  This time, by  using  the  method of papers \cite{Akhmet1,AkhmetBook,AkhmetExistence,Akhmet2020,Akhmet2024,iop,Akhmet2025} we investigate  existence  of the  new types of Poisson  stable solutions. 
 The newly invented method of verification of the Poisson stability   joined with the presence of the  periodic component  in the recurrence make possible  an extension  for the class  of studied differential equations.  In  papers \cite{Cheban1,Cheban2,Shcherbakov,Shcherbakov1} and others, quasilinear  systems   are  with  constant  matrices of coefficients,  and in our  case we research  systems with  periodic and,  even  with  Poisson  stable coefficients.   Another significant novelty,  which  is achieved in the present paper as well as in former our  studies \cite{Akhmet1,AkhmetBook,Akhmet2020,iop}  is the  numerical simulation of the Poisson  stable functions and solutions.  We believe that  altogether,  the present suggestions can  shape  a new  interesting science  direction,  not  only  in the  theoretical  study  of differential  equations,  but  also about  rich  opportunities for applications in  mechanics,  electronics,  artificial  neural  networks,  neuroscience. 
 
 Throughout the paper, $\mathbb R$ and $\mathbb{N}$ will stand for the set of real and natural numbers, respectively. Additionally,  the norm $\|u\|_1= \sup_{t\in \mathbb R} \|u(t)\|,$ where $\displaystyle \left\|u\right\|=\max_{1\le i \le n} \left|u_{i}\right|,$ $u=(u_{1},\ldots,u_{n}), u_{i} \in \mathbb{R}, i=1,2,...,n,$ will be used. Correspondly, for   a square matrix $A  = \{a_{ij}\},   i,j=1,2,...,n,$  the norm $\|A\|=\displaystyle \max_{i=1,\ldots,n} \displaystyle \sum_{j=1}^{n} |a_{ij}|$  will be utilized.

\textbf{Definition 1.} \label{poisson}\rm \cite{sell} A continuous and bounded function $\psi(t):\mathbb{R}\rightarrow\mathbb{R}^n$ is called Poisson stable, if there exists  a sequence $t_k,$  which diverges to infinity such that the sequence $\psi(t+t_k)$ converges to $\psi(t)$  uniformly on bounded intervals of $\mathbb{R}.$

 The sequence $t_k$ is said to be the \textit{Poisson sequence} of the function $\psi(t).$
 
 By Lemma \ref{tau} in the Appendix A, for a  positive fixed $\omega$ there exist a subsequence $t_{k_l}$ of the Poisson sequence $t_k$ and a  number $\tau_{\omega}$ such that $t_{k_l} \rightarrow \tau_{\omega} (mod \ \omega)$ as $l \rightarrow \infty.$ We shall call the number $\tau_{\omega}$ as the \textit{Poisson shift} for the Poisson sequence $t_k$ with respect to the $\omega.$ It is not difficult to find that for the fixed $\omega$ the set of all Poisson shifts, $T_{\omega},$ is not empty, and it can consist of several and even infinite number elements. The number $\kappa_{\omega}=inf\, T_{\omega},$  $0\leq \kappa_{\omega}<\omega,$ is said to be \textit{the Poisson number} for the Poisson sequence $t_k$ with respect to the number $\omega.$ 
 
\textbf{Definition 2.}\label{mppsf} \rm The sum $\phi(t)+\psi(t)$ is said to be a \textit{modulo periodic Poisson stable} ($MPPS$) function, if $\phi(t)$ is a continuous periodic and $\psi(t)$ is a Poisson stable functions.

 We shall  call  the function $\phi(t)$  the \textit{periodic component } and the function $\psi(t)$ the \textit{Poisson  component}   of the $MPPS$ function  in  what follows.

 \textbf{Remark 1.}  Duo to  Lemma \ref{mpps},  an $MPPS$ function  is a Poisson  stable if $\kappa_{\omega}$ equals zero.   Otherwise,    without  loss of generality,  the sequence  $\phi(t+t_k)+\psi(t+t_k)$  converges  on all  compact subsets of the real  axis to  the function  $\phi(t +\tau_{\omega})+\psi(t),$   where  $\tau_{\omega}$  is a nonzero Poisson shift  for the sequence $t_k.$  Since of the periodicity of the function $\phi(t),$ one can accept the last  convergence as a  special form  of recurrence. In  the next section, we shall  consider it   as a result  of  Theorem  \ref{linear}.

\section{Main results}

\subsection{Linear system of differential equations}

Consider the following system 
\begin{eqnarray} \label{system1}
	x'(t)=A(t)x(t)+\phi(t)+\psi(t),
\end{eqnarray}
where $t\in \mathbb{R},$ $x\in \mathbb{R}^n,$ $n\in \mathbb{N},$ $\phi(t):\mathbb{R}\rightarrow\mathbb{R}^n$ and $\psi(t):\mathbb{R}\rightarrow\mathbb{R}^n$ are continuous functions,  $A(t)$ is a continuous $n\times n$ matrix.

We assume that the following conditions are satisfied.
\begin{itemize}
	\item[\bf (C1)] $A(t)$ is an $\omega-$periodic matrix for a  fixed positive $\omega;$ 
	\item[\bf (C2)] $\phi(t)$ is an $\omega-$periodic function, and $\psi(t)$ is a Poisson stable function with a Poisson  sequence $t_k;$
	\item[\bf (C3)] the Poisson number $\kappa_{\omega}$ for the sequence $t_k$ is equal to zero.
\end{itemize}

According to Definition \ref{mppsf} and condition (C2), the sum $\phi(t)+\psi(t)$ is an $MPPS$ function. That is, the linear system (\ref{system1}) is with $MPPS$ perturbation.

Let us consider the homogeneous system, associated with (\ref{system1}),
\begin{eqnarray} \label{homo}
	x'(t)=A(t)x(t).
\end{eqnarray}
Let $X(t)$, $t \in \mathbb{R},$  is the fundamental matrix of the system (\ref{homo}) such that $X(0)=I,$  and $I$  is the $n\times n$ identical matrix.  Moreover,  $X(t,s)$ is transition matrix of the system (\ref{homo}),  which equal to $X(t)X^{-1}(s),$ and    $X(t+\omega,s+\omega)=X(t,s)$ for all $t,s \in \mathbb{R}.$

We assume that the following additional assumption is valid. 
\begin{itemize}
	\item[\bf (C4)] The multipliers of the system (\ref{homo}) are in modulus less than one.
\end{itemize}

It follows from the last condition that there exist positive numbers $ K\geq1$ and $\alpha$ such that
\begin{eqnarray}\label{exponent}
	\|X(t,s)\|\leq Ke^{-\alpha (t-s)}, 
\end{eqnarray}  
for $t\geq s$ \cite{Hartman}. 

\textbf{Lemma 1.}\label{estimation}  If the inequality (\ref{exponent}) is satisfied, then the following estimation is correct
	\begin{eqnarray}\label{transition}
		&&\|X(t+\tau,s+\tau)-X(t,s)\|\leq \displaystyle \max_{t \in \mathbb{R}}\|A(t+\tau)-A(t)\| \frac{2K^2}{\alpha^2e}e^{-\frac{\alpha}{2}(t-s)},
	\end{eqnarray}
	for $t\geq s$ and arbitrary real number $\tau.$  

\textbf{Proof.} Since
	\begin{eqnarray*}
		&&\displaystyle \frac{d X(t+\tau,s+\tau)}{d t}=A(t)X(t+\tau,s+\tau)+(A(t+\tau)-A(t))X(t+\tau,s+\tau),
	\end{eqnarray*}
	we have that
	\begin{eqnarray*}
		&&\displaystyle  X(t+\tau,s+\tau)=X(t,s)+\displaystyle \int_{s}^{t}X(t,u)(A(u+\tau)-A(u))X(u+\tau,s+\tau)du.
	\end{eqnarray*}
	That is why, 
	\begin{eqnarray*}
		&&\displaystyle  \|X(t+\tau,s+\tau)-X(t,s)\| \leq \\
		&&\displaystyle \int_{s}^{t}\|X(t,u)\|\|A(u+\tau)-A(u)\|\|X(u+\tau,s+\tau)\|du\leq \\
		&&\displaystyle \max_{t \in \mathbb{R}}\|A(t+\tau)-A(t)\| \int_{s}^{t}K^2e^{-\alpha(t-s)}du=\\
		&&\displaystyle \max_{t \in \mathbb{R}}\|A(t+\tau)-A(t)\|\frac{K^2}{\alpha}e^{-\alpha(t-s)}(t-s)= \\
		&&\displaystyle \max_{t \in \mathbb{R}}\|A(t+\tau)-A(t)\|\frac{K^2}{\alpha}e^{-\frac{\alpha}{2}(t-s)}e^{-\frac{\alpha}{2}(t-s)}(t-s).
	\end{eqnarray*}
	Since  $\displaystyle \sup_{u\geq 0}e^{-\frac{\alpha}{2}u}u=\frac{2}{\alpha e},$ the lemma is proved.

\textbf{Theorem 1.} \label{linear} Assume that conditions (C1), (C2) and (C4) are valid. Then the system (\ref{system1}) admits a unique asymptotically stable $MPPS$ solution.

\textbf{Proof.} The bounded solution of system (\ref{system1}) has the form \cite{Hartman}
	\begin{eqnarray}\label{integral}
		&&x(t)=\int_{-\infty}^t X(t,s)[\phi(s)+\psi(s)]ds,\ t \in \mathbb{R}.
	\end{eqnarray}
	
	One can write that $x(t)=x_{\phi}(t)+x_{\psi}(t),$ where $\displaystyle x_{\phi}(t)=\int_{-\infty}^t X(t,s)\phi(s)ds$ and $\displaystyle x_{\psi}(t)=\int_{-\infty}^t X(t,s)\psi(s)ds.$
	
	It is not difficult to show that the function $x_{\phi}(t)$ is $\omega-$periodic \cite{Farkas}.
	
	Next, we prove that the function $x_{\psi}(t)$ is Poisson stable. Fix arbitrary positive number $\epsilon$ and  interval $[a,b],$ $-\infty <a < b < \infty.$ We will show that for a large $k$ it is true that $\|x_{\psi}(t+t_k)-x_{ \psi}(t)\| <\epsilon$ on $[a,b].$ Let us choose two numbers $c$ and $\xi$ such that $c <a$ and $\xi$ is positive,  satisfying the following inequalities,
	\begin{eqnarray}\label{01}
		&&\displaystyle \frac{4K^2m_{\psi}}{\alpha^3 e} \xi <  \frac{\epsilon}{3},
	\end{eqnarray} 
	\begin{eqnarray}\label{11}
		&&\displaystyle \frac{2Km_{\psi}}{\alpha}e^{-\alpha(a-c)} <  \frac{\epsilon}{3},
	\end{eqnarray} 
	and
	\begin{eqnarray}\label{12}
		&&\displaystyle \frac{K\xi}{\alpha}  [1-e^{-\alpha(b-c)}] < \frac{\epsilon}{3},
	\end{eqnarray} 
	with $m_{\psi} = \displaystyle \sup_{t\in \mathbb R}\|\psi(t)\|.$ By applying condition (C4), without loss of generality,  for sufficiently  large $k$ we obtain that    $\|A(t+t_k)-A(t)\|<\xi$ for all $t\in \mathbb{R},$ and $\|\psi(t+t_k)-\psi(t)\| < \xi$ for $t\in[c,b].$ By using Lemma \ref{estimation} we attain that 
	\begin{eqnarray*}\label{Operatorr} 
		&& \|x_{\psi}(t+t_k)-x_{\psi}(t)\| 	= \|\int_{-\infty}^{t} \Big(X(t+t_k,s+t_k)\psi(s+t_k)-X(t,s)\psi(s)\Big)ds\|\leq  \\
		&& 	\int_{-\infty}^{t} \|X(t+t_k,s+t_k)-X(t,s)\|\|\psi(s+t_k)\|ds+ \\
		&&\int_{-\infty}^{t} \|X(t,s)\|\|\psi(s+t_k)-\psi(s)\|ds=  \\ 
		&&  \int_{-\infty}^{t} \|X(t+t_k,s+t_k)-X(t,s)\|\|\psi(s+t_k)\|ds+ \\
		&&\int_{-\infty}^{c} \|X(t,s)\|\|\psi(s+t_k)-\psi(s)\|ds+\int_{c}^{t} \|X(t,s)\|\|\psi(s+t_k)-\psi(s)\|ds \leq  \\
		&&\int_{-\infty}^{t} \frac{2K^2\xi}{\alpha^2 e}e^{-\frac{\alpha}{2}(t-s)}m_{\psi}ds+
		\int_{-\infty}^{t} 2Ke^{-\alpha(t-s)}m_{\psi}ds+
		\int_{-\infty}^{t} Ke^{-\alpha(t-s)}\xi ds\leq  \\
		&&	\frac{4K^2\xi}{\alpha^3 e} m_{\psi}+\frac{2Km_{\psi}}{\alpha}e^{-\alpha(a-c)}+\frac{K\xi}{\alpha} [1-e^{-\alpha(b-c)}].
	\end{eqnarray*}
	Now, the  inequalities (\ref{01}) to (\ref{12}) imply that $\|x_{\psi}(t+t_k)-x_{\psi}(t)\|<\epsilon,$ for $t\in [a,b].$ Therefore, the sequence $x_{\psi}(t+t_k)$ uniformly converges to $x_{\psi}(t)$ on each bounded interval. Thus, according to the Definition \ref{mppsf}  the solution $x(t)$ of the system (\ref{system1}) is $MPPS$ function with the periodic component   $x_{\phi}(t)$  and the Poisson component  $x_{\psi}(t).$ The asymptotic stability of the $MPPS$ solution   can be verified in the same way as  for the  bounded solution of a linear inhomogeneous system \cite{Farkas}.

\begin{figure}[ht] 
	\includegraphics[width=10.5 cm]{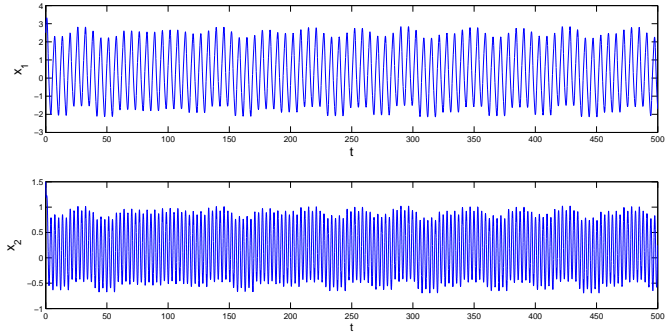} 	
	\caption{Coordinates of the  solution  $\bar{x}(t)$  of system (\ref{Ex1})  with initial values $\bar{x}_1(0)=2.5$ and $\bar{x}_2(0)=1.5,$  which converges  asymptotically to the $MPPS$ solution of system  (\ref{Ex1}).} 
	\label{fig1}
\end{figure}

\begin{figure}[ht] 	
	\includegraphics[width=10.5 cm]{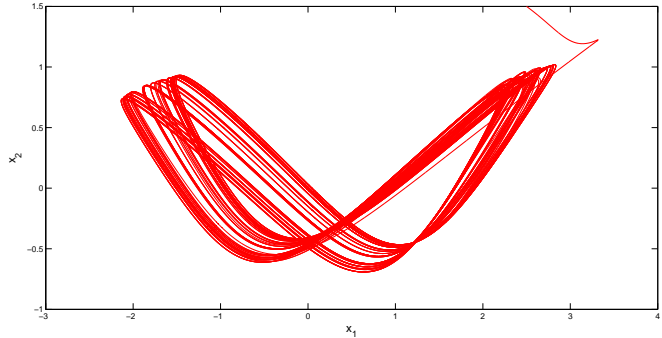} 	
	\caption{The trajectory of the solution $\bar{x}(t)$ of the equation (\ref{Ex1}).} 
	\label{fig2}
\end{figure}

\textit{Example 1.}  Let us consider the following linear inhomogeneous system, 
\begin{eqnarray} 
	\label{Ex1}
	\begin{array}{l}
		x_1'=(-1+0.5sin(2t))x_1+2.5cos(t)+5.5\Theta^2(t),\\
		x_2'=(-2+0.25cos(t))x_2+2sin(2t)+1.7\Theta(t),
	\end{array}
\end{eqnarray}
where $\Theta(t)=\int_{-\infty}^{t}e^{-3(t-s)}\Omega_{(3.85;6\pi)}(s)ds$ is the Poisson stable function described in Appendix B.   The perturbation is  an $MPPS$ function  with the periodic component  $\phi(t)= \Big(2.5cos(t),2sin(2t)\Big)^T$  and the Poisson component $\psi(t)= \Big(5.5\Theta^2(t),1.7\Theta(t)\Big)^T.$ The common period of the  coefficient  $A(t)$ and the periodic component   $\phi(t)$ is  $2\pi.$ Since the function $\Omega_{(3.85,6\pi)}(t)$ is constructed on the  intervals $[6\pi i,6\pi(i+1)),$ $i \in \mathbb{Z},$ for the Poisson sequence $t_k$ of the function $\Theta(t)$ there exists a subsequence $t_{k_l}$ such that $t_{k_l}\rightarrow 0 (mod \ 2\pi).$  Therefore, the Poisson number $\kappa_{\omega}=0.$ Condition (C4) is valid with the multipliers $\rho_1=e^{-2\pi},$ and $\rho_2=e^{-4\pi}.$   According to Theorem \ref{linear}, the system admits a unique asymptotically stable $MPPS$ solution. We depict in Figure \ref{fig1}  the coordinates of the solution  $\bar{x}(t),$  with initial values $\bar{x}_1(0)=2.5$ and $\bar{x}_2(0)=1.5,$ which 
asymptotically converges  to the $MPPS$ solution of the system  (\ref{Ex1}),  and consequently the coordinates visualize  the  solution.  In Figure \ref{fig2} the trajectory of the  solution $\bar{x}(t)$  is  shown.

\begin{figure}[ht] 	
	\includegraphics[width=10.5 cm]{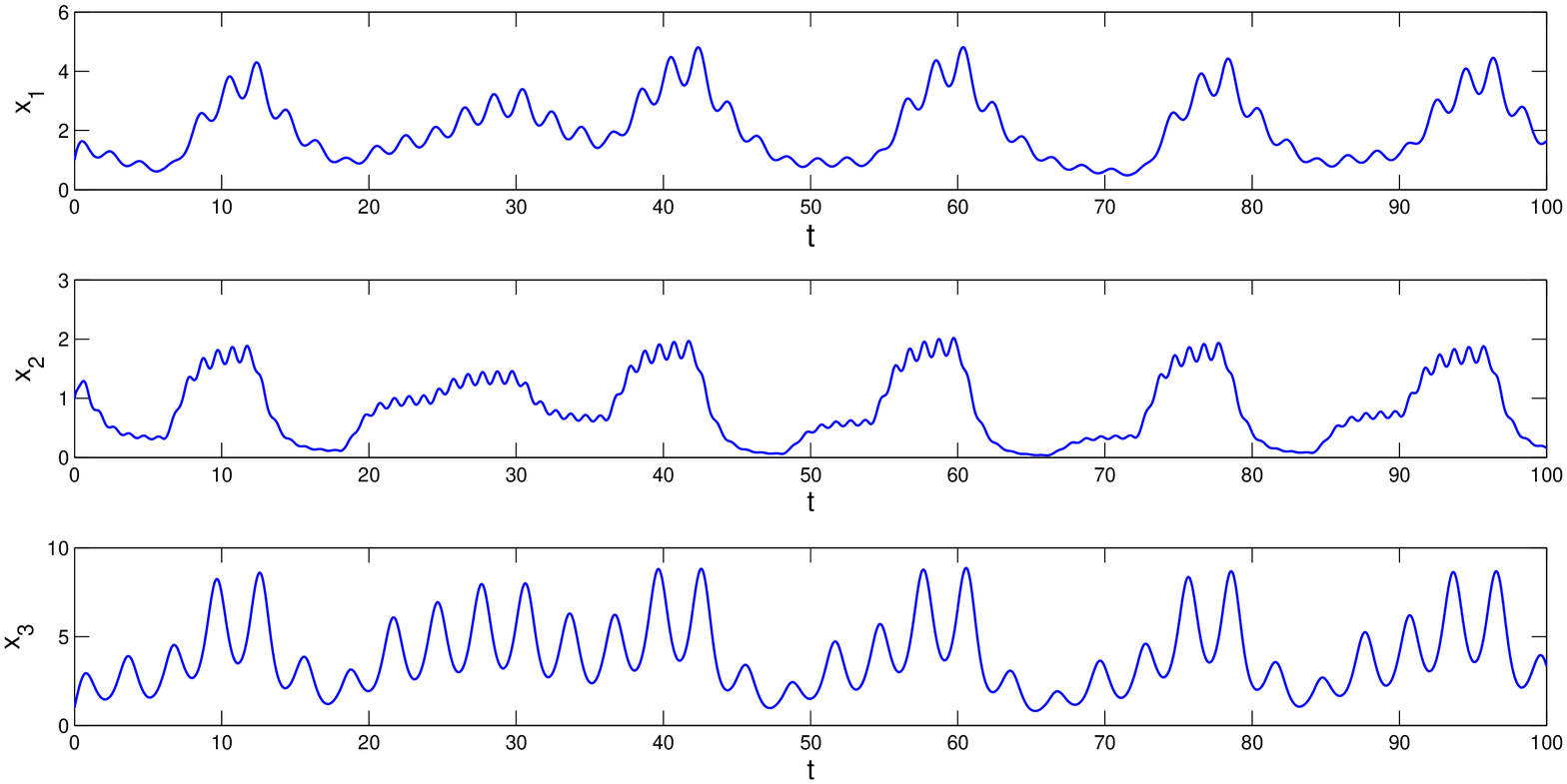} 
	\caption{Coordinates of the solution $\bar{x}(t),$ with initial values $\bar{x}_1(0)=1,$ $\bar{x}_2(0)=1$ and $\bar{x}_3(0)=1,$ which asymptotically converges to the $MPPS$ solution of the system (\ref{Ex2})}. 
	\label{fig3} 
\end{figure}

\begin{figure}[ht]
	\includegraphics[width=10.5 cm]{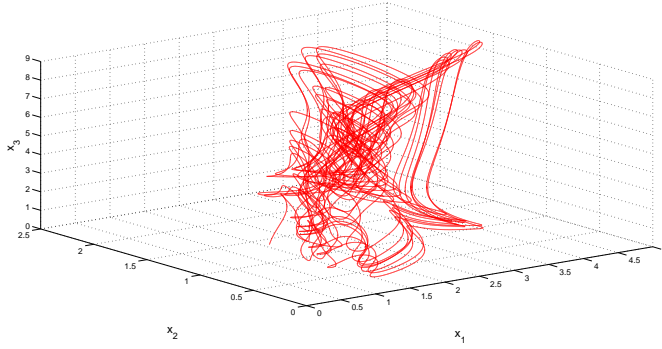}
	\caption{The trajectory of the  solution, $\bar{x}(t),$ of the equation (\ref{Ex2}).}
	\label{fig4}
\end{figure}

In the next example of the system (\ref{Ex1}),  the  periodic  component  $\phi(t)$ equals to the zero identically,  but  the condition $(C2)$ is correct,  since the  constant  function is  of  arbitrary period.  It is remarkable to  say  that the absence of a proper non-constant   periodic component  makes the dynamics more irregular,  what  is seen  in  Figures \ref{fig3} and \ref{fig4}. 

\textit{Example 2.}	Consider the inhomogeneous linear system  
\begin{eqnarray} 
	\label{Ex2}
	\begin{array}{l}
		x_1'=(-0.25+0.5cos(\pi t))x_1+12\Theta^3(t),\\
		x_2'=(-1.5+sin^2(\pi t))x_2+8\Theta^2(t),\\
		x_3'=(-0.5+cos(\frac{2\pi}{3}t))x_3+6\Theta(t),
	\end{array}
\end{eqnarray}
where $\Theta(t)=\int_{-\infty}^{t}e^{-2(t-s)}\Omega_{(3.9;6)}(s)ds.$ The matrix of coefficients is $3-$periodic, and conditions (C1)-(C3) are satisfied. The condition (C4) is valid with multipliers $\rho_1=e^{-0.75},$ $\rho_2=e^{-3}$ and $\rho_3=e^{-1.5}.$ Figure \ref{fig3} presents the coordinates of the  solution $\bar{x}(t)$ with initial values $\bar{x}_1(0)=1,$ $\bar{x}_2(0)=1$ and $\bar{x}_3(0)=1.$ The  solution $\bar{x}(t)$ approximates the coordinates of the $MPPS$ solution $x(t)$ of the equation (\ref{Ex2}) as time increases. The trajectory of the  solution $\bar{x}(t)$ is shown in Figure \ref{fig4}.

\subsection{Quasilinear differential equations}

The main object of the present section is the system of quasilinear differential equations
\begin{eqnarray}\label{quasi}
	x'(t)=A(t)x+g(t,x)+\phi(t)+\psi(t),
\end{eqnarray}
where $t \in \mathbb{R}, x \in \mathbb R^n,$ $n$ is a fixed natural number;  $A(t)$ is $n-$dimensional square matrix and satisfies to the condition (C1) and inequality (\ref{exponent});  $g:\mathbb R\times U \rightarrow \mathbb{R}^n, g =(g_1,\ldots,g_n),$ $U=\{x \in \mathbb{R}^n, \|x\|<H\},$ where $H$ is a fixed positive number; the functions $\phi(t)$ and $\psi(t)$ satisfy conditions (C2) and (C3). 

The following conditions on system (\ref{quasi}) are required.

\begin{itemize}
	\item[\bf (C5)] the function $g(t,x)$ is continuous and $\omega-$periodic in $t;$
	\item[\bf (C6)] there exists a positive constant $L$ such that $\left\|g(t,x_1)-g(t,x_2)\right\| 
	\le L \left\| x_1-x_2 \right\|$  for   all $t \in \mathbb{R}, x_1, x_2 \in U.$
\end{itemize}

We denote $\displaystyle \sup_{\mathbb{R}\times U} \|g(t,x)\|=m_g,$ $\displaystyle \max_{t \in \mathbb{R}} \|\phi(t)\|=m_{\phi}$ and $\displaystyle \sup_{t \in\mathbb{R}} \|\psi(t)\|=m_{\psi}.$

The following additional conditions will be needed:
\begin{itemize}
	\item[\bf (C7)] $\displaystyle  \frac{K(m_g+m_{\phi}+m_{\psi})}{H}<\alpha;$
	\item[\bf (C8)] $KL<\alpha.$
\end{itemize}
For simplicity, we use the notation $F(t,x)=g(t,x)+\phi(t)+\psi(t) $ in what follows.

According to \cite{Hartman}, a bounded on the real axis function $y(t)$ is a solution of (\ref{quasi}),  if and only if it satisfies the equation
\begin{eqnarray}\label{dif_eqn}
	&&y(t)=\int_{-\infty}^t X(t,s)F(s,y(s))ds,\ t \in \mathbb{R}.
\end{eqnarray}

\textbf{Theorem 2.}\label{quasilinear} If conditions (C1)-(C8) are valid, then the system (\ref{quasi}) possesses a unique asymptotically stable Poisson stable solution. 

\textbf{Proof.} Let $t_k$ is the Poisson sequence of the function $\psi(t)$ in the system (\ref{quasi}). We denote $B$ the set of all Poisson stable  functions $\nu(t)=(\nu_1, \nu_2,..., \nu_n),$ $\nu_i \in\mathbb R,$ $i=1,2,...,n,$ with common Poisson sequence $t_k,$ which satisfy $\left\|\nu \right\|_{1} < H.$  
	
	Let us show that the   $ B$ is a complete space.  Consider  a Cauchy  sequence $\theta_m(t)$ in $ B$,   which   converges to  a limit  function $\theta(t)$ on $\mathbb{R}$.      We have that
	\begin{eqnarray}\label{Poisson}
		&&\|\theta(t+t_k)-\theta(t)\| < \|\theta(t+t_k)-\theta_m(t+t_k)\|+\|\theta_m(t+t_k)-\theta_m(t)\|+\nonumber \\
		&&\|\theta_m(t)-\theta(t)\|. 	
	\end{eqnarray}
	for a fixed closed and bounded interval $I \subset \mathbb{R}.$ Now, one can take sufficiently large $m$ and $k$ such that each term on the right hand-side of (\ref{Poisson}) is smaller than $\frac{\epsilon}{3}$ for a fixed positive $\epsilon$ and $t\in I$. That is, the sequence $\theta(t+t_k)$ uniformly converges to $\theta(t)$  on $I.$
	Likewise, one   can    check that  the limit    function  is uniformly  continuous \cite{Hartman}. The completeness of $ B$ is shown.
	
	Define the operator $\Pi$ on $ B$ such that 
	\begin{eqnarray}\label{oper}
		&&\Pi \nu(t)=\int_{-\infty}^t X(t,s)F(s,\nu(s))ds, \ t \in \mathbb{R}.
	\end{eqnarray}
	Fix a function $\nu(t)$ that belongs to $ B$. We have that 
	\begin{eqnarray*}
		&& \|\Pi \nu(t)\| \leq  \displaystyle \int_{-\infty}^t  \|X(t,s)\| \|F(s,\nu(s))\| ds 
		\leq \displaystyle \frac{K(m_g+m_{\phi}+m_{\psi})}{\alpha}
	\end{eqnarray*}
	for all $t\in\mathbb R.$ 
	Therefore, by the condition (C7) it is true that $\left\|\Pi \nu \right\|_1 <H$.
	
	Fix a positive number $\epsilon$ and an interval $[a,b],$ $-\infty <a < b < \infty.$ Let us choose two numbers $c<a,$ and $\xi>0$ satisfying the inequalities
	\begin{eqnarray}\label{q1}
		\frac{4K^2\xi}{\alpha^3 e} (m_g+m_{\phi}+m_{\psi})<  \frac{\epsilon}{3},	
	\end{eqnarray}
	\begin{eqnarray}\label{q2}
		&&\displaystyle \frac{2K}{\alpha}(m_g+m_{\phi}+m_{\psi})e^{-\alpha(a-c)} <  \frac{\epsilon}{3},
	\end{eqnarray} 
	and 
	\begin{eqnarray}\label{q3}
		&&\displaystyle \frac{K\xi}{\alpha}  [1-e^{-\alpha(b-c)}] < \frac{\epsilon}{3}.
	\end{eqnarray}
	By using condition (C4) and Lemmas \ref{mpps},  \ref{poissonstable} in the Appendix, without loss of generality,  we obtain that  $\|A(t+t_k)-A(t)\|<\xi$ for all $t \in \mathbb{R},$ and $\|F(t+t_k,\nu(t+t_k)) - F(t,\nu(t))\| < \xi$  for $t\in[c,b]$ and  sufficiently  large $k.$   Then, applying the inequality  (\ref{transition}), we get that  
	\begin{eqnarray*}\label{Operator} 
		&& \|\Pi \nu(t+t_k) -\Pi \nu(t)\|	= \\
		&&\|\int_{-\infty}^{t} X(t+t_k,s+t_k)F(s+t_k,\nu(s+t_k))ds-\int_{-\infty}^{t}X(t,s)F(s,\nu(s))ds\|	\leq  \\
		&&\int_{-\infty}^{t} \|X(t+t_k,s+t_k)-X(t,s)\|\|F(s+t_k,\nu(s+t_k))\|ds+ \\
		&&\int_{-\infty}^{c}\|X(t,s)\|\| F(s+t_k,\nu(s+t_k))-F(t,s)\|ds \leq  \\ 
		&&  \int_{c}^{t}\|X(t,s)\|\| F(s+t_k,\nu(s+t_k))-F(t,s)\|ds \leq  \\
		&&\int_{-\infty}^{t} \frac{2K^2\xi}{\alpha^2 e}e^{-\frac{\alpha}{2}(t-s)}(m_g+m_{\phi}+m_{\psi})ds+\\
		&&
		\int_{-\infty}^{t} 2Ke^{-\alpha(t-s)}(m_g+m_{\phi}+m_{\psi})ds+
		\int_{-\infty}^{t} Ke^{-\alpha(t-s)}\xi ds\leq  \\
		&&\frac{4K\xi}{\alpha^3 e} (m_g+m_{\phi}+m_{\psi})+ \frac{2K}{\alpha}(m_g+m_{\phi}+m_{\psi}) e^{-\alpha(a-c)} + \frac{K\xi}{\alpha} [1-e^{-\alpha(b-c)}],
	\end{eqnarray*}
	for all  $t \in [a,b].$ From inequalities (\ref{q1}) to (\ref{q3}) it follows that
	$
	\|\Pi \nu(t+t_k) -\Pi \nu(t)\| < \epsilon
	$
	for $t\in[a,b].$ Therefore, $\Pi \nu(t+t_k)$ uniformly converges to $\Pi \nu(t) $  on bounded interval of $\mathbb R.$ 
	
	It is easy to  verify  that $\Pi \nu(t)$ is a  uniformly continuous function, since its derivative is a uniformly bounded function on the real axis. Summarizing the above discussion,  the set $ B$ is invariant for the operator $\Pi$.
	
	We proceed to  show that the operator $\Pi : B\rightarrow  B$ is contractive. Let $u(t)$ and $v(t)$ be members of $ B$. Then, we obtain  that
	\begin{eqnarray*}
		&& \|\Pi u(t) -\Pi v(t)\| 
		\leq \displaystyle \int_{-\infty}^t \|X(t,s)\|\|F(s,u(s)) -F(s,v(s))\| ds\leq \\
		&& \displaystyle \int_{-\infty}^t Ke^{\alpha(t-s)}L \|u(s)-v(s)\|ds \leq \displaystyle \frac{KL}{\alpha} \|u(t)-v(t)\|_1
	\end{eqnarray*}
	for all $t \in \mathbb R.$ Therefore,  the inequality $\left\|\Pi u - \Pi v\right\|_1 \leq \displaystyle \frac{KL}{\alpha} \left\|u-v\right\|_1$ holds, and according to the condition $(C8)$ the operator $\Pi:  B \to  B$ is contractive.
	
	By contraction mapping theorem there exists  the unique fixed point, $x(t) \in B,$  of the operator $\Pi,$ which is the unique bounded Poisson stable solution of the system (\ref{quasi}).

	Finally, we will study the asymptotic stability of the Poisson solution $x(t)$ of the system (\ref{quasi}). It is true that 
	\begin{eqnarray*}
		x(t)=X(t,t_0)x(t_0)+\int_{t_0}^{t}X(t,s)\Big(g(s,x(s))+\phi(s)+\psi(s)\Big)ds,
	\end{eqnarray*}
	for $t\geq t_0.$
	
	Let $z(t)$ be another solution of system (\ref{quasi}). One can write
	\begin{eqnarray*}
		z(t)=X(t,t_0)z(t_0)+\int_{t_0}^{t}X(t,s)\Big(g(s,z(s))+\phi(s)+\psi(s)\Big)ds.
	\end{eqnarray*}
	Making use the relation
	\begin{eqnarray*}
		x(t)-z(t)=X(t,t_0)(x(t_0)-z(t_0))+\int_{t_0}^{t}X(t,s)\Big(g(s,x(s))-g(s,z(s))\Big)ds,
	\end{eqnarray*}
	we obtain that
	\begin{eqnarray*}
		&&\|x(t)-z(t)\|\leq\|X(t,t_0)\|\|x(t_0)-z(t_0)\|+\int_{t_0}^{t}\|X(t,s)\|\|g(s,x(s))-g(s,z(s)\|ds\leq \\
		&&Ke^{-\alpha(t-t_0)}\|x(t_0)-z(t_0)\|+\int_{t_0}^{t}KLe^{-\alpha(t-s)}\|x(s)-z(s)\|ds.
	\end{eqnarray*}
	Now, applying Gronwall-Bellman Lemma, one can attain that
	\begin{eqnarray}\label{stable}
		&&\|x(t)-z(t)\|\leq Ke^{-(\alpha-KL)(t-t_0)}\|x(t_0)-z(t_0)\|.
	\end{eqnarray}
	The last inequality and condition (C8) confirm that the Poisson stable solution $x(t)$ is  asymptotically stable. The theorem is proved. $\Box$

\textbf{Remark 2.} According to the Lemma \ref{psmpps} in the Appendix, the Poisson stable solution $x(t)$ of the system (\ref{quasi}) is an $MPPS$ function.

\textit{Example 3.} Consider the quasilinear system.
\begin{eqnarray} 
	\label{Ex3}
	&&	\begin{array}{l}
		x_1'=(-1.5+2sin(2t))x_1+0.01cos(2t)arctg(x_2)+1.2sin(8t)-10.5\Theta^3(t),\\
		x_2'=(-3.5+3sin^2(2t))x_2+0.03sin(4t)arctg(x_3)-1.5cos(8t)+2.5\Theta(t),\\
		x_3'=(-1.5+2cos^2(t))x_3-0.02sin(2t)arctg(x_1)+sin(4t)+7.2\Theta^2(t),
	\end{array}
\end{eqnarray}
where $\Theta(t)=\int_{-\infty}^{t}e^{-3(t-s)}\Omega_{(3.86,3\pi)}(s)ds$ is the Poisson stable function, which described similarly to  that  in  Appendix B. Since, the piecewise constant function $\Omega_{(3.86;3\pi)}(t)$ is given on intervals $[3\pi i, 3\pi(i+1)),$ for the Poisson sequence $t_k$ of the function $\Theta(t)$ there exists a subsequence $t_{k_l}$ such that $t_{k_l}\rightarrow 0(mod \ \pi),$ that is the condition (C3) is valid. The common period of the matrix $A(t)$ and functions $g(t,x),$ $\phi(t)$ is equal to $\pi.$  We have that the function $g(t,x)=(0.01cos(2t)arctg(x_2), 0.03sin(4t)arctg(x_3), -0.02sin(2t)arctg(x_1))^T$ is continuous and $\pi-$ periodic in $t$ and satisfies  condition (C6) with $L=0.03.$  The sum of  $\phi(t)=(1.2sin(8t), -1.5cos(8t), sin(4t))^T$ and $\psi(t)=(10.5\Theta^3(t), 2.5\Theta(t), 7.2\Theta^2(t))^T$ is an $MPPS$ function, which meets conditions (C2),(C3).    The assumptions (C4)-(C8) are valid with $m_g=0.048,$ $m_{\phi}=1.5,$ $m_{\psi}=0.84,$ $\rho_1=e^{-1.5\pi},$ $\rho_2=e^{-2\pi},$ $\rho_3=e^{-0.5\pi},$ $\alpha=0.5\pi,$ $K=1,$ and $H=4.8.$   Thus,  all conditions for the last  theorem have been  verified,   and there is the Poisson stable solution of the system,  which is asymptotically stable.  

\begin{figure}[ht] \includegraphics[width=10.5cm]{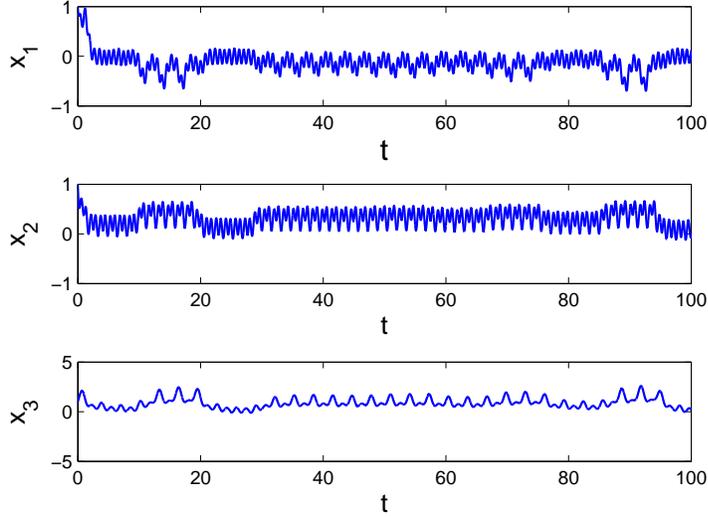} 	\caption{The coordinates of the  solution $\bar{x}(t),$ which    are  asymptotic  for   the Poisson stable solution of the system (\ref{Ex3}). } 
	\label{fig5}
\end{figure}

\begin{figure}[ht]
	\includegraphics[width=10.5cm]{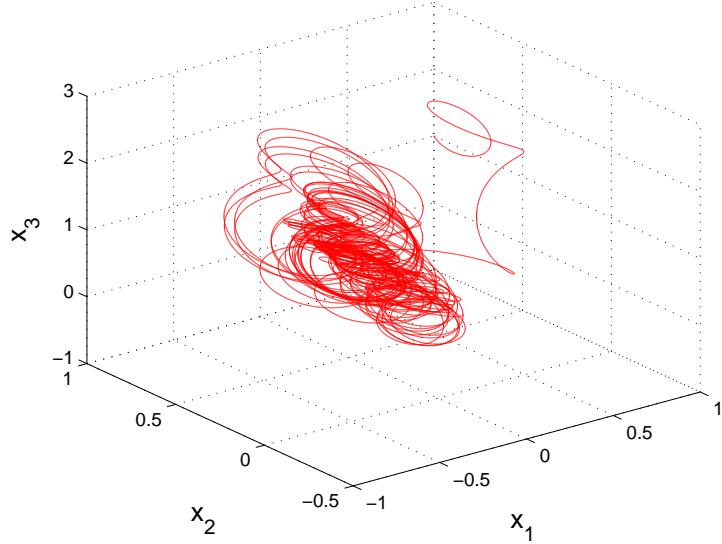}
	\caption{The trajectory of the  solution  $\bar{x}(t),$  which  illustrates  the Poisson stability of the system (\ref{Ex3}).}
	\label{fig6}
\end{figure}

It is worth noting that the simulation of the Poisson stable solution,  $x(t),$  is not possible, since the initial value is not known precisely. For this reason, we will consider the  solution $\bar{x}(t)$ of the system (\ref{Ex3}), with initial values $\bar{x}_1(0)=1,$ $\bar{x}_2(0)=1$ and $\bar{x}_3(0)=1.$ Using the inequality (\ref{stable}) one can obtain that $\|x(t)-\bar{x}(t)\|\leq e^{-1.54}\|x(0)-\bar{x}(0)\|$ for $t\geq 0.$ The last inequality shows that  $\|x(t)-\bar{x}(t)\|$ decreases exponentially. Consequently, the graph of the  solution $\bar{x}(t)$  asymptotically approaches to the Poisson stable solution $x(t)$ of the system (\ref{Ex3}), as time increases. The Figure \ref{fig5} demonstrates the coordinates of the solution  $\bar{x}(t),$  which illustrate the Poisson stability of the system  (\ref{Ex3}).   In the Figure \ref{fig6} the trajectory of the function $\bar{x}(t)$ is depicted. 

\section{A case with MPPS coefficients}

Let us consider the quasilinear equation (\ref{quasi}) with $A(t)=B(t)+D(t),$ where $B(t)$ is a continuous $\omega-$periodic matrix, and $D(t)$ is a Poisson stable matrix with the Poisson sequence $t_k.$ That is, the coefficient is an  $MPPS$ matrix and  the system (\ref{quasi}) is of the  form
\begin{eqnarray}\label{quasi2}
	x'(t)=(B(t)+D(t))x+g(t,x)+\phi(t)+\psi(t),
\end{eqnarray}
where the functions $\phi(t)$ and $\psi(t)$ satisfy conditions (C2) and (C3) and their sum is an  $MPPS$ function. The function $g(t,x)$ satisfies conditions (C5), (C6).

Denote $G(t,x)=D(t)x+g(t,x)+\phi(t)+\psi(t)$ and rewrite the system (\ref{quasi2}) as
\begin{eqnarray}
	x'(t)=B(t)x+G(t,x).
\end{eqnarray}
The homogeneous  $\omega-$ periodic system, associated  with (\ref{quasi2}),
\begin{eqnarray}\label{homoquasi}
	y'(t)=B(t)y,
\end{eqnarray}
has the  fundamental matrix $Y(t),$ $Y(0)=I,$ and the transition matrix $Y(t,s),$ $t,s \in \mathbb{R}.$

Assume that the following assumptions are valid.
\begin{enumerate}
	\item[\bf (C9)] The multipliers of the system (\ref{homoquasi}) are in modulus less than one.
\end{enumerate}
From the condition (C9) we have that there exist positive numbers $D\geq 1$ and $\beta$ such that
\begin{eqnarray}
	\|Y(t,s)\|\leq De^{-\beta(t-s)},
\end{eqnarray}
for $t\geq s.$

\begin{itemize}
	\item[\bf (C10)]\quad $D(L+d)<\beta;$ 
	\item[\bf (C11)]\quad $\displaystyle  \frac{D(m_g+m_{\phi}+m_{\psi})}{H}<\beta-Dd,$
\end{itemize}
where $d=\sup_{t \in \mathbb{R}}\|D(t)\|.$

The next theorem is proved in similar for theorem  \ref{quasilinear} way.

\textbf{Theorem 3.} If conditions (C2), (C3), (C5), (C6),  and (C9)  to (C11) are hold, then system (\ref{quasi2}) admits a unique asymptotically stable Poisson stable solution. 

\section{Appendix A}

\textbf{Lemma 2.} \label{tau} For arbitrary sequence of positive real numbers $t_k,$ $k=1,2,\cdots,$ and a positive number $\omega$ there exist a subsequence $t_{k_l},$ $l=1,2,\cdots,$ and a  number $\tau_{\omega},$ $0\leq \tau_{\omega}<\omega,$ such that $t_{k_l} \rightarrow \tau_{\omega} (mod \ \omega)$ as $l \rightarrow \infty.$

\textbf{Proof.} Consider the sequence $\tau_k$ such that $t_k\equiv \tau_k(mod \ \omega),$ and $0\leq\tau_k<\omega$ for all $k\geq 1.$ The boundedness of the sequence $\tau_k$ implies that there exist a  subsequence $\tau_{k_l},$  which converges to a number $\tau_{\omega}$ \cite{Haggarty}.

\textbf{Lemma 3.} \label{kappa} $\kappa_{\omega} \in T_{\omega}.$

\textbf{Proof.} Assume on the contrary that $\kappa_{\omega}$ is not in $ T_{\omega}.$ Then there exists a strictly decreasing sequence $\tau_m,$  $m\geq 1,$ in $ T_{\omega},$ such that $\tau_m \rightarrow \kappa_{\omega}.$ For each natural $m,$ denote by $t^m_i$  a subsequence of $t_k$ such that $t^m_i \rightarrow \tau_m (mod \ \omega)$ as $i \rightarrow \infty.$
	
	Fix a sequence of positive numbers $\epsilon_n,$ which converges to the zero. One can find numbers $i_n,$ $n=1,2,\ldots,$ such that 
	$ |t^{n}_{i_{n}}-\tau_{n}|<\epsilon_{n} (mod \ \omega).$ It is clear that $t^{n}_{i_{n}}\rightarrow \kappa_{\omega} (mod \ \omega)$ as $n \rightarrow \infty.$

\textbf{Remark 3.}  The last  assertion implies that if  $\kappa_{\omega}=0,$ then there exists a subsequence $t_{k_l}$ such that $t_{k_l} \rightarrow 0 (mod \ \omega)$ as $l\rightarrow \infty.$

\textbf{Lemma 4.} \label{mpps} If $f(t)=\phi(t)+\psi(t)$ is an $MPPS$  function, and $\kappa_{\omega}=0,$ then the function $f(t)$ is a Poisson stable.

\textbf{Proof.} According to Lemma \ref{kappa}, there exists a subsequence $t_{k_l,}$ which tends to zero in modulus $\omega$ as $l \rightarrow \infty.$ Without loss of generality assume that $t_k \rightarrow 0(mod \ \omega)$ as $k\rightarrow \infty.$ Fix a positive number $\epsilon,$ and bounded interval $I\subset \mathbb{R}.$ The periodic function $\phi(t)$ is uniformly continuous on $\mathbb{R}.$ Consequently, there exists a number $k_1$ such that 
	\begin{eqnarray*}
		\|\phi(t+t_{k})-\phi(t)\| < \frac{\epsilon}{2}, 
	\end{eqnarray*}
	for all $t \in \mathbb{R}$ and $k>k_1.$ Moreover, there exists an integer $k_2,$ such that 
	\begin{eqnarray*}
		\|\psi(t+t_{k})-\psi(t)\|< \frac{\epsilon}{2},
	\end{eqnarray*}
	for $t \in I,$ $k>k_2.$ This is why, 
	\begin{eqnarray*}
		&&\|f(t+t_{k})-f(t)\|\leq\|\phi(t+t_{k})-\phi(t)\|+\|\psi(t+t_{k})-\psi(t)\|<\epsilon,
	\end{eqnarray*}
	if $t \in I$ and $k>\max(k_1,k_2).$
	That is, the function $f(t)$ is Poisson stable. 

\textbf{Lemma 5.}\label{psmpps}
	Assume that  $\psi(t)$ is a Poisson stable function. If $\kappa_{\omega}=0,$ for some positive number $\omega,$   then $\psi(t)$ is  an $MPPS$ function.

\textbf{Proof.} Let us write $\psi(t)=g(t)+(\psi(t)-g(t)),$ where $g(t)$ is a continuous $\omega-$periodic function. Since $\kappa_{\omega}=0,$ then the subtraction $\psi(t)-g(t)$	is Poisson stable by Lemma \ref{mpps}.

\textbf{Remark 4.} The last result  is a source for  the optimization problem how to choose  the function $g(t)$  and the period $\omega$ to  minimize the difference  $\psi(t)-g(t).$  In other words,  the problem  of  approximation of Poisson stable functions  with periodic ones. 	It  is of exceptional interest  for celestial mechanics \cite{Poin}.

\textbf{Lemma 6.} \label{poissonstable} Assume that a function $G(t,u):\mathbb{R}\times U \rightarrow \mathbb{R}^n, U\subseteq \mathbb{R}^n,$ is a Poisson stable function  in $t$ and satisfies the inequality $\left\|G(t,u_1)-G(t,u_2)\right\| 
	\le L \left\| u_1-u_2 \right\|,$ where $L$ is a positive constant,  for   all $t \in \mathbb{R}, u_1, u_2 \in U.$ Moreover, $\upsilon(t):\mathbb{R} \rightarrow U$ is  $\omega-$periodic in $t.$  If the Poisson sequence and period $\omega$ are  such that the Poisson number $\kappa_{\omega}$ equals to the  zero, then the function $G(t,\upsilon(t))$ is Poisson stable.

\textbf{Proof.} By the Lemma \ref{kappa}  there exists a subsequence $t_{k_l},$ such that $t_{k_l} \rightarrow 0 (mod \ \omega)$  as $l \rightarrow \infty.$ We assume, without loss of generality, that the sequence $t_k$ itself satisfies the condition $t_{k} \rightarrow 0 (mod \ \omega)$  as $k \rightarrow \infty.$ 
	
	Let us fix a positive number $\epsilon,$ and a bounded interval $I.$ Since $t_{k} \rightarrow 0 (mod \ \omega)$  as $k \rightarrow \infty,$ for sufficiently large  $k,$ we obtain that $\displaystyle \|G(t+t_{k},\upsilon(t+t_{k}))-G(t,\upsilon(t+t_{k}))\|<\frac{\epsilon}{2}$ for all $t\in \mathbb{R},$ and $\displaystyle \|\upsilon(t+t_{k})-\upsilon(t)\|<\frac{\epsilon}{2L}$ for $t\in I.$ We have that
	\begin{eqnarray*}	&&\|G(t+t_{k},\upsilon(t+t_{k}))-G(t,\upsilon(t))\|\leq \|G(t+t_{k},\upsilon(t+t_{k}))-G(t,\upsilon(t+t_{k}))\|+\\	&&\|G(t,\upsilon(t+t_{k}))-G(t,\upsilon(t))\|\leq \frac{\epsilon}{2}+L\frac{\epsilon}{2L}\leq \epsilon,\end{eqnarray*}
	for all $t\in I.$ That is, $G(t,\upsilon(t))$ is the Poisson stable function. 

Similarly one can prove the following assertions.

\textbf{Lemma 7.}\label{spoissonstable} Assume that a function $G(t,u):\mathbb{R}\times U \rightarrow \mathbb{R}^n, U\subseteq \mathbb{R}^n,$ is $\omega-$periodic in $t$ and satisfies the inequality $\left\|G(t,u_1)-G(t,u_2)\right\| 
	\le L \left\| u_1-u_2 \right\|,$ where $L$ is a positive constant,  for   all $t \in \mathbb{R}, u_1, u_2 \in U.$ Moreover, $\upsilon(t):\mathbb{R} \rightarrow U$ is a Poisson stable function.  If the Poisson sequence  and period $\omega$ are   such that the Poisson number $\kappa_{\omega}$ equals  to  the zero, then the function $G(t,\upsilon(t))$ is Poisson stable.

\textbf{Lemma 8.}\label{ppoissonstable} Assume that a function $G(t,u):\mathbb{R}\times U \rightarrow \mathbb{R}^n, U\subseteq \mathbb{R}^n,$ is  Poisson stable  in $t$  and satisfies the inequality $\left\|G(t,u_1)-G(t,u_2)\right\| 
	\le L \left\| u_1-u_2 \right\|,$ where $L$ is a positive constant,  for   all $t \in \mathbb{R}, u_1, u_2 \in U.$ Moreover, $\upsilon(t):\mathbb{R} \rightarrow U$ is a Poisson stable function. If there exists  a  Poisson sequence  common  for the functions $G(t,u)$ and $\upsilon(t),$  then the function $G(t,\upsilon(t))$ is Poisson stable.

\textbf{Remark 5.} The last  lemma implies,  in particular,   that  sum and  product   of Poisson stable functions with common  Poisson sequence are Poisson  stable functions.

\section{Appendix B}

This part  of the paper is about   an  example of the Poisson stable  functions.  The  task is not  easy one,  and there  very  few   constructively  determined  cases \cite{sell,Nemytskii}.  In our  research,  we utilize the dynamical approach  of functions determination.  One of the  most  familiar  is of sin  and cos functions as solutions of ordinary  differential equations.  We  shall consider the Poisson  function as a continuous component  of solution for a hybrid  system,  which  consists of   a discrete equation and a simple differential equation,  while  discrete  component  can  be accepted  as a Poisson  stable sequence.  A significant  element   of the present  study  is visualization of the continuous  Poisson stable solution through  a neighboring it  by an asymptotically  close    counterpart. 

In \cite{ak8}  as a part  of the result   construction  of  a Poisson stable  sequence  was performed as the  solution of the logistic  equation
\begin{eqnarray} \label{log}
	&&\lambda_{n+1}= \mu \lambda_{n}(1-\lambda_{n}) .
\end{eqnarray}

More precisely, it is proved that for each $\mu\in[3+(2/3)^{1/2},4]$  there exists a solution $\{\eta_n\},$ $n \in \mathbb{Z},$ of equation (\ref{log}) such that the sequence belongs to the interval $[0,1]$  and there exists a sequence $\zeta_n,$ which diverges to infinity such that $|\eta_{i+\zeta_n}-\eta_{i}|\rightarrow 0$ as $n \rightarrow \infty$ for each $i$ in bounded intervals of integers.

Consider the following integral
\begin{eqnarray}
	\Theta(t)=\int_{-\infty}^{t}e^{-2(t-s)}\Omega(s)ds, \ t\in \mathbb{R},
\end{eqnarray}
where $\Omega(t)$ is a piecewise constant function defined on the real axis through the equation $\Omega(t)=\eta_i$ for $t \in [i,i+1),$ $i\in \mathbb{Z}.$ It  is  convenient to  consider the function $\Theta(t)$  as a unique bounded on the real  axis  solution  of the  equation $\Theta^{\prime} = -2 \Theta + \Omega(t).$  In  all next  examples  of the paper we use the function notation $\Omega(t) =  \Omega_{(p,q)}(t),$ where $p$ is value of constant $\mu,$ and $q$ denotes the length of the intervals on which the function $\Omega(t)$ is built.

It is worth noting that $\Theta(t)$ is bounded on the hole real axis such that $\sup_{t \in \mathbb{R}}|\Theta(t)|\leq 1/2.$

Next, we will show that $\Theta(t)$ is a Poisson stable function.

Consider a fixed closed interval $[a,b]$ of the axis and a positive number $\varepsilon.$ Without loss of generality one can assume that $a$ and $b$ are integers. Let us fix a positive number $\xi$ and an integer $c<a,$ which satisfy the following inequalities $e^{-2(a-c)}<\frac{\varepsilon}{2}$ and $\xi[1-e^{-2(b-c)}]<\varepsilon.$ Let $n$ be a large natural number such that $|\Omega_{(3.89,1)}(t+\zeta_n)-\Omega_{(3.89,1)}(t)|<\xi$ on $[c,b].$ Then for all $t \in[a,b]$ we obtain that
\begin{eqnarray*}
	&&|\Theta(t+\zeta_n)-\Theta(t)|=|\int_{-\infty}^{t}e^{-2(t-s)}(\Omega_{(3.89,1)}(s+\zeta_n)-\Omega_{(3.89,1)}(s))ds|= \\
	&&|\int_{-\infty}^{c}e^{-2(t-s)}(\Omega_{(3.89,1)}(s+\zeta_n)-\Omega_{(3.89,1)}(s))ds+\\
	&&\int_{c}^{t}e^{-2(t-s)}(\Omega_{(3.89,1)}(s+\zeta_n)-\Omega_{(3.89,1)}(s))ds|\leq\\
	&&\int_{-\infty}^{c}e^{-2(t-s)}2ds+\int_{c}^{b}e^{-2(t-s)}\xi ds\leq e^{-2(a-c)}+\frac{\xi}{2}[1-e^{-2(b-c)}]<\frac{\varepsilon}{2}+\frac{\varepsilon}{2}=\varepsilon.
\end{eqnarray*}
Thus, $|\Theta(t+\zeta_n)-\Theta(t)|\rightarrow 0$ as $n\rightarrow \infty$ uniformly on the interval $[a,b].$

\section{Conclusions}  In  this paper, we have introduced a new type of recurrence,  which  is a   sum  of the two components,   the rigorously  periodic,  which  is an  oscillatory by the definition,     and  the Poisson  stable as a recurrent  motion.  Accordingly, we call   it as modulo periodic Poisson stable function. Sufficient condition for the  dynamics   to  be Poisson stable has been  determined.    The novelty is  convenient  for theoretical  analysis  as solutions  of  differential equations of different types, as well  discrete equations.  This time,  we restrict ourself  with  pioneering results  for quasi-linear  ordinary  differential  equations.    If  one remove the periodic component  from the discussions of the present  paper,   
the results  became identical  to  those obtained  previously \cite{Akhmet1,AkhmetBook,AkhmetExistence,Akhmet2020,Akhmet2024,iop,Akhmet2025}.  If one take in  consideration the periodicity  in the Poisson stability,  and those achievements of the paper  for simulations of the recurrence,  the results   open  new productive opportunities  in  the  research  of mechanical,  electronic dynamics, and we mostly  rely in  neuroscience.   Concerning  theoretical research,  it  is of big  interest  to  search  for Poisson  stability  and its periodic components  in  such  famous dynamics as Lorenz,  R\"{o}ssler and Chua attractors.  Generally speaking,  one can  search  for periodic  components  of  any chaotic dynamics.   This definitely can  be joined with the problems of optimization, if one looks for minimal  contribution of the   recurrent component.

\textbf{Acknowledgments} M. Akhmet and A. Zhamanshin have  been supported by  2247-A National Leading Researchers Program of TUBITAK, Turkey, N 120C138.
	M. Tleubergenova  has been supported by the Science Committee of the Ministry of Education and Science of the Republic of Kazakhstan (grant No. AP08955400 and No. AP08856170).

\end{document}